\magnification=\magstep1
\baselineskip =5.5mm
\lineskiplimit =1.0mm
\lineskip =1.0mm

\long\def\comment#1{}

\def\writemonth#1{\ifcase#1
\or January\or February\or March\or April\or May\or June\or July%
\or August\or September\or October\or November\or December\fi}

\newcount\mins
\newcount\minmodhour
\newcount\hour
\newcount\hourinmin
\newcount\ampm
\newcount\ampminhour
\newcount\hourmodampm
\def\writetime#1{%
\mins=#1%
\hour=\mins \divide\hour by 60
\hourinmin=\hour \multiply\hourinmin by -60
\minmodhour=\mins \advance\minmodhour by \hourinmin
\ampm=\hour \divide\ampm by 12
\ampminhour=\ampm \multiply\ampminhour by -12
\hourmodampm=\hour \advance\hourmodampm by \ampminhour
\ifnum\hourmodampm=0 12\else \number\hourmodampm\fi
:\ifnum\minmodhour<10 0\number\minmodhour\else \number\minmodhour\fi
\ifodd\ampm p.m.\else a.m.\fi
}

\font\sevenrm=cmr7
\font\fiverm=cmr5
\font\sevenex=cmex10 at 7pt
\font\sevenbf=cmbx7
\font\fivebf=cmbx5
\font\sevenit=cmti7
\font\seventeenrm=cmr12 at 17pt
\font\twelverm=cmr12
\font\seventeeni=cmmi12 at 17pt
\font\twelvei=cmmi12
\font\seventeensy=cmsy10 at 17pt
\font\twelvesy=cmsy10 at 12pt
\font\seventeenex=cmex10 at 17pt
\font\seventeenbf=cmbx12 at 17pt
\def\sevensize{\sevenrm \baselineskip=4.5 mm%
\textfont0=\sevenrm \scriptfont0=\fiverm \scriptscriptfont0=\fiverm%
\def\rm{\fam0 \sevenrm}%
\textfont1=\seveni \scriptfont1=\fivei \scriptscriptfont1=\fivei%
\def\mit{\fam1 } \def\oldstyle{\fam1 \seveni}%
\textfont2=\sevensy \scriptfont2=\fivesy \scriptscriptfont2=\fivesy%
\def\cal{\fam2 }%
\textfont3=\sevenex \scriptfont3=\sevenex \scriptscriptfont3=\sevenex%
\def\bf{\fam\bffam\sevenbf} \textfont\bffam\sevenbf
\scriptfont\bffam=\fivebf \scriptscriptfont\bffam=\fivebf
\def\it{\fam\itfam\sevenit} \textfont\itfam\sevenit
}
\def\seventeensize{\seventeenrm \baselineskip=5.5mm%
\textfont0=\seventeenrm \scriptfont0=\twelverm \scriptscriptfont0=\sevenrm%
\def\rm{\fam0 \seventeenrm}%
\textfont1=\seventeeni \scriptfont1=\twelvei \scriptscriptfont1=\seveni%
\def\mit{\fam1 } \def\oldstyle{\fam1 \seventeeni}%
\textfont2=\seventeensy \scriptfont2=\twelvesy \scriptscriptfont2=\sevensy%
\def\cal{\fam2 }%
\textfont3=\seventeenex \scriptfont3=\seventeenex%
\scriptscriptfont3=\seventeenex%
\def\bf{\fam\bffam\seventeenbf} \textfont\bffam\seventeenbf
}

\def\setheadline #1, #2 \par{\headline={\ifnum\pageno=1 
\hfil
\else \sevensize \noindent
\ifodd\pageno \hfil #2\hfil \else
\hfil #1\hfil \fi\fi}}

\footline={\ifnum\pageno=1 \fiverm \hfil
This paper is not in final form. Typeset using \TeX\ on
\writemonth\month\ \number\day, \number\year\ at \writetime{\time}\hfil 
\else \rm \hfil \folio \hfil \fi}

\def\finalversion{\footline={\ifnum\pageno=1 \hfil 
\else \rm \hfil \folio \hfil \fi}}

\def\Head #1:{\medskip\noindent{\bf #1:}}
\def\Proof:{\medskip\noindent{\bf Proof:}}
\def\Proofof #1:{\medskip\noindent{\bf Proof of #1:}}
\def\endproof{\nobreak\hfill$\sqr$\bigskip\goodbreak}
\def\itemi{\item{i)}}
\def\itemii{\item{ii)}}
\def\itemiii{\item{iii)}}

\def\dt{\it}

\def\ts#1{{\textstyle{#1}}}
\def\Abstract\par#1\par{\centerline{\vtop{
\sevensize
\abovedisplayskip=6pt plus 3pt minus 3pt
\belowdisplayskip=6pt plus 3pt minus 3pt
\moreabstract\parindent=0 true in% 
A{\fiverm BSTRACT}: \ \ #1}}
\abovedisplayskip=12pt plus 3pt minus 9pt
\belowdisplayskip=12pt plus 3pt minus 9pt
\vskip 0.4 true in}
\def\moreabstract{%
\par \hsize = 5 true in \hangindent=0 true in \parindent=0.5 true in}
\outer\def\firstbeginsection#1\par{\bigskip\vskip\parskip
\message{#1}\leftline{\bf#1}\nobreak\smallskip\noindent}

\def\sqr{\vcenter {\hrule height.3mm
\hbox {\vrule width.3mm height 2mm \kern2mm
\vrule width.3mm } \hrule height.3mm }}

\def\Bbb{\bf}
\def\E{{\Bbb E}}
\def\R{{\Bbb R}}
\def\Z{{\Bbb Z}}
\def\N{{\Bbb N}}
\def\C{{\Bbb C}}
\font\eightrm=cmr10 at 8pt
\def\R{\hbox{\rm I\kern-2pt R}}
\def\Z{\hbox{\rm Z\kern-3pt Z}}
\def\N{\hbox{\rm I\kern-2pt I\kern-3.1pt N}}
\def\C{\hbox{\rm \kern0.7pt\raise0.8pt\hbox{\eightrm I}\kern-4.2pt C}} 
\def\E{\hbox{\rm I\kern-2pt E}}

\def\invp{{1\over p}}
\def\invq{{1\over q}}
\def\invc{{c^{-1}}}
\def\half{{1\over 2}}
\def\smallhalf{\ts\half}

\def\list#1,#2{#1_1$, $#1_2,\ldots,$\ $#1_{#2}}
\def\lists#1{#1_1$, $#1_2,\ldots}

\def\set#1{\{1$, $2,\ldots,$\ $#1\}}

\def\lnorm{\left\|}
\def\rnorm{\right\|}
\def\normo#1{\lnorm #1 \rnorm}

\def\lmod{\left|}
\def\rmod{\right|}
\def\modo#1{\lmod #1 \rmod}

\def\angleo#1{\left\langle #1 \right\rangle}

% document relevent defs

\def\F{{\cal F}}
\def\cN{{\cal N}}

\def\invN{{1\over N}}
\def\kot{{k\over2}}

\def\TliNY{T:\liN\to Y}
\def\TCKY{T:C(K)\to Y}
\def\TXY{T:X\to Y}

\def\pitoT{\pi_{2,1}(T)}
\def\pipqT{\pi_{p,q}(T)}
\def\pitT{\pi_2(T)}
\def\pipT{\pi_p(T)}
\def\betatT{\beta^{(2)}(T)}

\def\KtT{K^{(2)}(T)}

\def\bx{\bar{x}}

\def\xS{\list x,S}
\def\epsilonlist{\lists\varepsilon}

\def\bxS{\list{\bx},S}
\def\yS{\list y,S}
\def\zS{\list z,S}

\def\liN{{l_\infty^N}}

\def\Lpq{{L_{p,q}}}

\def\sumS{\sum_{s=1}^S}
\def\sumN{\sum_{n=1}^N}

\def\randwalk#1{\normo{\sumS #1_s x_s}}
\def\Bernwalk{\randwalk\varepsilon}

\def\Bernwalki{\Bernwalk_\infty}

\def\rootmean#1,#2,#3,#4{\left(\sumS\normo{#1}_{#2}^{#3}\right)^#4}

\def\pthmpTxs{\rootmean T(x_s),,p,\invp}

\def\wkqnormxs{\sup \left\{ \left( \sumS \modo{\angleo{x^*,x_s}}^q
               \right)^\invq \right\} }

% document starts here

\setheadline THE RADEMACHER COTYPE OF OPERATORS FROM $l_\infty^N$,
             MONTGOMERY-SMITH --- TALAGRAND

\finalversion

{
\seventeensize
\centerline{\bf The Rademacher Cotype of}
\centerline{\bf Operators from $\liN$}
}
\bigskip\bigskip
\centerline{\bf S.J.~Montgomery-Smith}
{
\sevensize\baselineskip=4.0mm
\centerline{\it Department of Mathematics, University of Missouri,}
\centerline{\it Columbia, MO 65211.}
}
\medskip
\centerline{\bf M.~Talagrand}
{\sevensize\baselineskip=4.0mm
\centerline{\it Department of Mathematics, The Ohio State University,}
\centerline{\it 231 W.\ 18th Avenue, Columbus, OH 43210.}
\vskip 1pt
\centerline{\it Equipe d'Analyse -- Tour 46, Universit\'e Paris VI,}
\centerline{\it 4 Place Jussieu, 75230 Paris Cedex 05.}
}
\bigskip

\Abstract

We show that for any operator $\TliNY$, where $Y$\ is a Banach
space, that its cotype~2 constant, $\KtT$, is related to its
$(2,1)$-summing
norm, $\pitoT$, by
$$ \KtT \le c \, \log\log N \,\, \pitoT .$$
Thus, we can show that there is an operator $\TCKY$\ that has cotype~2, but is
not 2-summing.
\medskip\moreabstract\noindent{\bf A.M.S.\ Classification:} Primary 46B20,
Secondary 60G99

\firstbeginsection Introduction

The notation we use in this paper is loosely based on that given in 
[L--T1], [L--T2] and [P1].

We let $\epsilonlist$\ be independent
Rademacher random variables, that is, $\Pr(\varepsilon_s=1) =
\Pr(\varepsilon_s=-1) = \smallhalf $.
A linear operator $\TXY$\ is said to have {\dt (Rademacher) cotype~$p$\/}
($p\ge2$) if there is a constant $C<\infty$\ such that for all $\xS$\
in $X$\ we
have 
$$ \pthmpTxs \le C \,\, \E \Bernwalk. $$
The smallest value of $C$\ is called the {\dt (Rademacher)
cotype~$p$\ constant\/} of $T$, and is denoted by $K^{(p)}(T)$. These
definitions extend to spaces in the obvious way; a space $X$\ has
cotype~$p$\
if its identity operator has cotype~$p$.

We define the {\dt $(p,q)$-summing norm\/} of a linear operator
$\TXY$, denoted
by $\pipqT$, to be the least number $C$\ such that for all $\xS$ in
$X$\ we have 
$$ \pthmpTxs \le C \, \wkqnormxs, $$ 
where the supremum is taken over all $x^*$\ in the unit ball of $X^*$.
We call a $(p,p)$-summing operator a {\dt $p$-summing operator\/},
and write $\pipT$\ for $\pi_{p,p}(T)$. We say that the operator is
{\dt $(p,q)$-summing\/} ({\dt $p$-summing\/}) if $\pipqT < \infty$
(respectively $\pipT < \infty$).

If $1\le p<\infty$, and $1\le q \le \infty$, then we let $\Lpq(\mu)$\
denote the
Lorentz space on the measure $\mu$. We
refer the reader to [H] or [L--T2] for details, but just note that the
$L_{p,1}$\ norm may be calculated using
$$ \normo{f}_{p,1} = \int_0^\infty \mu(\modo f > t)^\invp \,dt 
   = \ts{\invp} \int_0^\infty s^{\invp-1} f^*(s) \,ds ,$$
where $f^*$\ denotes the non-decreasing rearrangement of $\modo f$.

The basic motivation behind this paper is in classifying operators
from $C(K)$\
that factor through a Hilbert space, where $C(K)$\ denotes the
continuous
functions on the compact Hausdorff topological space, $K$. The first result in
this direction is due to Grothendieck,
which
states that any bounded linear operator $C(K) \to L_1$\  factors
through Hilbert
space. This was generalized by Maurey [Ma1], allowing $L_1$\ to be
replaced by any space of cotype~$2$, to give the following result (see also
[P1]).

\proclaim Theorem 1. Let $\TCKY$\ be a linear operator, where $Y$\ is
any Banach
space. Then the following are equivalent:
\item{i)} $T$\ is $2$-summing;
\item{ii)} $T$\ factors through Hilbert space;
\item{iii)} $T$\ factors through a space of cotype~$2$.

However, we are still left with the following question:
if the {\it operator\/} $\TCKY$\ has cotype~$2$, does it follow that
it factors
through Hilbert space?

One way one might tackle this problem is to consider the
$(2,1)$-summing norms
of such operators. Jameson [J] showed that there is an operator
$\TliNY$\ such
that $\pitT \ge \invc \sqrt{\log N} \, \pitoT$. Hence, if we can
establish a
strong relationship between the cotype~$2$\ constants and the
$(2,1)$-summing
norms of such operators, then we can answer the above question in the
negative.
To this end, we have the following --- the main result of this paper.

\proclaim Theorem 2. There is a constant $c$\ such that for any
operator
$\TliNY$, where $Y$\ is a Banach space, then the
cotype~$2$\ constant is bounded according to the relation:
$$ \KtT \le c \, \log\log N \, \pitoT .$$

\proclaim Corollary. There is an operator $\TCKY$, where $Y$\ is a
Banach space,
that has cotype~$2$, but does not factor through Hilbert space.

Finally, before embarking on the proof of this result, we point out
that for
$p>2$, the above problems have been completely answered.

\proclaim Theorem 3. Let $\TCKY$\ be a bounded linear operator, where
$Y$\ is a
Banach space. Then for all $p>2$, the following are equivalent:
\itemi $T$\ is $(p,1)$-summing;
\itemii $T$\ has Rademacher cotype~$p$;
\itemiii $T$\ factors through a space with Rademacher cotype~$p$.

The implication (i) $\Leftrightarrow$\ (ii) is due to Maurey [Ma2]. The third
equivalence follows from the fact that any $(p,1)$-summing operator from $C(K)$\
factors through $L_{p,1}$ (see [P2] or Theorem~5 below), and that $L_{p,1}$\ has
Rademacher cotype~$p$, (see  [C]).

\proclaim Theorem 4. If $p>2$, then there is a bounded linear
operator $C(K) \to
L_p$\ that is not $p$-summing.

We refer the reader to [K].

\beginsection Proof of the Main Result

To prove Theorem~2, we need the following two results. The first
allows us to
reduce questions about $(p,1)$-summing operators from $C(K)$\ to the
canonical
embedding $C(K) \to L_{2,1}(K,\mu)$\ ($\mu$\ a probability measure),
and is due
to Pisier (see [P2]).

\proclaim Theorem 5.  Let $\TCKY$\ be a $(p,1)$-summing operator,
where $Y$\ is
a Banach space, and $p \ge 1$. Then there is a Radon probability
measure
$\mu$\ on $K$\ and a constant $C \le p^\invp \pi_{p,1}(T) $\ such
that for all
$x\in C(K)$\ we have $\normo{Tx} \le C \, \normo x_{L_{p,1}(K,\mu)} $.

The second result is about Rademacher processes, and is due to the
second named
author (for the proof, see [Ld--T]). First we establish some more
notation. If
$T$\ is a bounded subset of $\R^S$, we write
$$ r(T) = \E \, \sup_{t\in T} \modo{\sumS \varepsilon_s t(s)} . $$
If $B$\ is a subset of $\R^S$, we write $\cN(T,B)$\ for the minimal
number
of translates of $B$\ required to cover $D$.
We write $B_1^S$\ for the unit ball of $l^S_1$, and
$B_2^S$\ for the unit ball of $l_2^S$. From now on, we take all
logarithms to
base $2$.

\proclaim Theorem 6. There is a constant $c_1$\ such that if $T$\ is
a bounded
subset of $\R^S$, and $\epsilon>0$, then letting $D=c_1\,r(T)\,B_1^S +
\epsilon \, B_2^S$, we have
$$ r(T) \ge c_1^{-1} \epsilon \sqrt{\log\cN(T,D)} .$$

Now we will state the main result towards proving Theorem~2.

\proclaim Proposition 7. There is a constant $c_2$\ such that if
$(\Omega,\F,\mu)$\ is a probability space with $N$\ atoms, and $\xS \in
L_\infty(\mu)$\ are such that  
$$ \E \Bernwalki \le 1 ,$$ 
then
$$ \rootmean{x_s},{L_{2,1}(\mu)},2,\half \le c_2\, \log\log N .$$

Our first step in establishing this result is to restate Theorem~6 in
a more
suitable form.

\proclaim Lemma 8. There is a constant $c_1$\ (the same one as in
Theorem~6) such that the following holds. Suppose that
$(\Omega,\F,\mu)$\ is a
measure space, with $\Omega$\ finite, and $\xS\in L_\infty(\mu)$\ with
$$ \E\Bernwalki \le 1 .$$
Then for all integers $k$, we may partition $\Omega$\ into
at most $2^{2^k}$\ measurable sets, find $\yS$, $\zS \in
L_\infty(\mu)$, and find
$\bxS\in L_\infty(\Omega,\F',\mu)$\ (where $\F'$\ denotes the algebra
generated
by the partition), such that $x_s = \bx_s + y_s + z_s$,
$$ \E\normo{\sumS \varepsilon_s \bx_s}_\infty \le 1
   , \quad
   \normo{\sumS\modo{y_s}}_\infty \le c_1
   \hbox{\quad and\quad}
   \normo{\left(\sumS\modo{z_s}^2\right)^\half}_\infty \le c_1 \,
2^{-\kot} .$$

\Proof: Let $T = \left\{\,\bigl(x_s(\omega)\bigr)_{s=1}^S : \omega
\in \Omega
\right\}$, and let $\epsilon = c_1^{-1} 2^{-\kot}$. If we apply
Theorem~6, we
see that there are $2^{2^k}$\ translates, $t_l + c_1(B_1^S + 2^{-\kot}
B_2^S)$\ ($1\le l \le 2^{2^k}$), that cover $T$.
We let the covering of $\Omega$\ be the sets 
$$ A_l = \left\{\, \omega : \bigl(x_s(\omega)\bigr)_{s=1}^S \in t_l + 
         c_1(B_1^S + 2^{-\kot} B_2^S) \right\} ,$$
and if $A_l$\ is non-empty, we choose $\omega_l \in A_l$. Define
$ \bx_s(\omega) = x_s(\omega_l)$\ if $\omega \in A_l$. Now, if
$\omega\in A_l$,
we know that $\bigl(x_s(\omega)-\bx_s(\omega)\bigr)_{s=1}^S \in
c_1(B_1^S +
2^{-\kot} B_2^S)$, that is, there are $\bigl(y_s(\omega)\bigr)_{s=1}^S
\in
c_1\,B_1^S$\ and $\bigl(z_s(\omega)\bigr)_{s=1}^S \in c_1\,2^{-\kot}
B_2^S$, with
$x_s(\omega) = \bx_s(\omega) + y_s(\omega) + z_s(\omega)$.
\endproof

\proclaim Lemma 9. There is a constant $c_3$\ such that if
$(\Omega,\F,\mu)$\
is a measure space with $\Omega$\ finite, then the following hold. 
\item{i)} If $y\in L_\infty(\mu)$, then $\normo
y_{2,1} \le \normo y_\infty^\half \normo y_1^\half $.
\item{ii)} If the smallest atom is of size $a$, then for all $z\in
L_\infty(\mu)$\ we have 
$$ \normo z_{2,1}\le c_3\, \left(1+\sqrt{\log(\mu(\Omega)/a)}\right) \normo z_2
.$$  \item{iii)} If there are $N$\ atoms, then for all $z\in
L_\infty(\mu)$\ we have $\normo z_{2,1}\le \sqrt N \normo z_2$.

\Proof: i) We have that
$$ \eqalignno{
   \normo{y}_{2,1}
   &= \int_0^{\normo{y}_\infty} \sqrt{\mu(\modo{y}>t)} \,dt \cr
   &\le \left(\int_0^{\normo{y}_\infty} \,dt\right)^\half
        \left(\int_0^{\normo{y}_\infty} \mu(\modo{y}>t)
\,dt\right)^\half
        \cr
   &= \normo{y}_\infty^\half \normo{y}_1^\half .\cr}$$

\noindent ii) We have
$$ \eqalignno{
   \normo{z}_{2,1}
   &= \ts{\half} \int_0^\infty {z^*(s)\over\sqrt s} \,ds
        \cr
   &\le \sqrt a \normo{z}_\infty +
        \ts{\half} \int_{a}^{\mu(\Omega)} {z^*(s)\over\sqrt s} \,ds
\cr
   &\le \normo{z}_2 +
        \ts{\half} \left(\int_{a}^{\mu(\Omega)}  \,{ds \over s}
\right)^\half
        \left(\int_{a}^{\mu(\Omega)} \bigl(z^*(s)\bigr)^2 \,ds
        \right)^\half \cr
   &\le c_3 \,\left(1+\sqrt{\log (\mu(\Omega)/a)}\right) \normo{z}_2 .\cr}$$

\noindent iii) Let $\list B,N$\ be the atoms of $\Omega$
arranged so that $z^*(n)$, the value of $\modo{z}$\ on $B_n$, is in
non-increasing order. Also, let $z^*(N+1) = 0 $. Then
$$ \eqalignno{
   \normo{z}_{2,1}
   &= \sumN \left(\sum_{m=1}^n \mu(B_m)\right)^\half
      \bigl(z^*(n)-z^*(n+1)\bigr) \cr
   &\le \sqrt N \left(\sumN \sum_{m=1}^n \mu(B_m)
        \bigl(z^*(n)-z^*(n+1)\bigr)^2 \right)^\half \cr
   &\le \sqrt N \left(\sumN \mu(B_n)
        \bigl(z^*(n)\bigr)^2 \right)^\half \cr
   &= \sqrt N \normo{z}_2 .\cr}$$
as desired.
\endproof

We remark that Lemma~9(i) is a well known interpolation result, and is true for
all measure spaces.

\proclaim Lemma 10. If $(\Omega,\F,\mu)$\ is a probability space with
$\Omega$\
finite, then 
$$ \rootmean{y_s},{2,1},2,\half \le \normo{\sumS\modo{y_s}}_\infty .$$

\Proof: This follows straight away from Lemma~9(i).
\endproof

Lemma~10 is also well known (and true for all probability spaces). In fact it is
a reformulation of the statement that the canonical embedding $C(\Omega) \to
L_{2,1}(\mu)$\ has $(2,1)$-summing norm equal to $1$.

\proclaim Lemma 11. There is a constant $c_4$\ such that, if
$(\Omega,\F,\mu)$\
is a probability space with at most $N$\ atoms, then 
$$ \rootmean{z_s},{2,1},2,\half \le c_4 \sqrt{\log N}
   \normo{\left(\sumS\modo{z_s}^2\right)^\half}_\infty .$$

\Proof: Let $A\subset\Omega$\ be the union of those atoms of measure
less than
${1\over N^2}$, so that $\mu(A) \le \invN$. By Lemma~9(ii), we have
that $\normo{z_s \chi_{\Omega\setminus A}}_{2,1} \le c_3 \sqrt{\log N}
\normo{z_s}_2$, and by Lemma~9(iii), we have that $\normo{z_s
\chi_A}_{2,1} \le
\sqrt N \normo{z_s \chi_A}_2 $. Thus, we have that
$$ \eqalignno{
   \rootmean{z_s},{2,1},2,\half 
   &\le \rootmean{z_s \chi_{\Omega\setminus A}},{2,1},2,\half 
        + \rootmean{z_s \chi_A},{2,1},2,\half  \cr
   &\le c_3\sqrt{\log N} \rootmean{z_s},2,2,\half 
        + \sqrt {N \, \mu(A)} 
        \normo{\left(\sumS\modo{z_s}^2\right)^\half}_\infty \cr
   &\le c_4\sqrt{\log N}
\normo{\left(\sumS\modo{z_s}^2\right)^\half}_\infty 
        ,\cr}$$
as desired.
\endproof

\Proofof Proposition~7: Without loss of generality, we may suppose
that
$N=2^{2^{k}}$. We prove the result by induction over $k$. Suppose that
$\Omega$\ has $2^{2^{k+1}}$\ atoms. Apply Lemma~8 to cover $\Omega$\
by
$2^{2^k}$\ subsets, and to give $\bxS$, $\yS$, $\zS$\ as described in
the lemma.
Then, by the triangle inequality 
$$ \rootmean{x_s},{2,1},2,\half 
   \le \rootmean{\bx_s},{2,1},2,\half 
   + \rootmean{y_s},{2,1},2,\half 
   + \rootmean{z_s},{2,1},2,\half .$$
By the induction hypothesis,
$$ \rootmean{\bx_s},{2,1},2,\half \le c_2 k .$$
By Lemmas~10 and~11 we have that
$$ \rootmean{y_s},{2,1},2,\half \le c_1 $$
and
$$ \rootmean{z_s},{2,1},2,\half \le c_1 c_4 \, 2^{-\kot} \left(1+\sqrt{\log
   2^{2^k}}\right)
   \le 2 c_1 c_4 .$$
Hence
$$ \rootmean{x_s},{2,1},2,\half \le c_2(k+1) ,$$
as required, taking $c_2 = 1 + 2 c_1 c_4$.
\endproof

To prove the main result is now easy.

\Proofof Theorem~2: By Theorem~5, it is sufficient to show that for
any
probability measure $\mu$\ on $\set N$, the cotype~$2$\ constant of
the
canonical embedding $\liN \to L_{2,1}(\mu)$\ is bounded by some
universal
constant times $\log\log N$. But this is precisely what Proposition~7
says.
\endproof

\beginsection Final Remarks

There is a similar result for Gaussian cotype (see [Mo2]).

\proclaim Theorem 12. There is a constant $c$\ such that, for any
operator
$\TliNY$, where $Y$\ is a Banach space, the
Gaussian cotype~$2$\ constant, $\betatT$, is bounded according to the
relation:
$$ \betatT \le c \sqrt{\log\log N} \, \pitoT .$$

This result is the best possible, as is implicitly shown in [T].

\proclaim Theorem 13. There is a constant $c$\ such that for any
integer $N$,
there is an operator $\TliNY$, where $Y$\ is a Banach space, such that
$$ \betatT \ge \invc \sqrt{\log\log N} \, \pitoT .$$

Since the Rademacher cotype~2 constant is greater than a constant
times the
Gaussian cotype~2 constant, we have the following corollary.

\proclaim Corollary. There is a constant $c$\ such that for any
integer $N$,
there is an operator $\TliNY$, where $Y$\ is a Banach space, such that
$$ \KtT \ge \invc \sqrt{\log\log N} \, \pitoT .$$

We also have the following, the result originally stated in [T].

\proclaim Corollary. There is an operator $\TCKY$, where $Y$\ is a
Banach
space, that is $(2,1)$-summing, but does not have Rademacher cotype~$2$.

If we write $R_N$\ for the supremum of $\KtT / \pitoT$\ over all
$\TliNY$, then
we have shown that $\invc \sqrt{\log\log N} \le R_N \le c \, \log\log
N$.
Clearly, we are left with the following problem.

\proclaim Open Question. What is the asymptotic behavior of $R_N$?

\beginsection Acknowledgements

The main result of this paper originally appears in the Ph.D.\ thesis
of the
first named author [Mo1], and he would like to express his
thanks to his advisor, D.J.H.~Garling, and to the Science and
Engineering
Council who financed his studies. He would also like to express his
gratitude to G.J.O.~Jameson who first suggested the problem to him,
and gave him
much encouragement.

\beginsection References

\halign{#\hfil & \quad\vtop{\hsize=5.5 true in%
\parindent=0pt\hangindent=1em
\strut#\strut} \cr
C & J.~Creekmore,\rm\ Type and cotype in Lorentz $L_{p,q}$\
spaces,\sl\
Indag.\ Math.\ {\bf 43} (1981), 145--152.\cr
H & R.A.~Hunt,\rm\ On $L(p,q)$\ spaces,\sl\ L'Enseignement Math.\ (2)
{\bf 12} (1966), 249--275.\cr
J & G.J.O.~Jameson,\rm\ Relations between summing norms of mappings
on
$l_\infty^n$,\sl\ Math.\ Z.\ {\bf 194} (1987), 89--94.\cr
K & S.~Kwapien,\rm\ On a theorem of L.~Schwartz and its applications
to
absolutely summing operators,\sl\ Stud.\ Math.\ {\bf 38} (1984),
193--200.\cr
Ld--T & M.~Ledoux and M.~Talagrand,\sl\ Isoperimetry and Processes in
Probability in a Banach Space,\rm\ Springer-Verlag (to appear).\cr
L--T1 & J.~Lindenstrauss and L.~Tzafriri,\sl\ Classical Banach Spaces
I---Se\-qu\-ence
Spa\-ces,\rm\ Springer-Verlag, 1977.\cr
L--T2 & J.~Lindenstrauss and L.~Tzafriri,\sl\ Classical Banach Spaces
II---Fu\-nc\-t\-ion
Spa\-ces,\rm\ Springer-Verlag, 1979.\cr
Ma1 & B.~Maurey,\sl\ Th\'eor\`emes de factorisation pour les op\'erateurs
lin\'eaires \`a valeurs dans un espace $L^p$,\rm\ Ast\'erisque~{\bf 11},
1974.\cr
Ma2 & B.~Maurey,\rm\ Type et cotype dans les espaces munis de
structures
locales
inconditionelles,\sl\ Seminaire Maurey-Schwartz 1973--74,
Exp.~24--25.\cr
Mo1 & S.J.~Montgomery-Smith,\sl\ The Cotype of Operators from
$C(K)$,\rm\
Ph.D.\ thesis, Cambridge, August 1988.\cr
Mo2 & S.J.~Montgomery-Smith,\rm\ The Gaussian cotype of operators from
$C(K)$,\sl\ Israel J.\ of Math.\ {\bf 68} (1989), 123--128.\cr
P1 & G.~Pisier,\sl\ Factorization of Linear Operators and Geometry
of Banach
Spaces,\rm\ Amer.\ Math.\ Soc., 1986.\cr
P2 & G.~Pisier,\rm\ Factorization of operators through $L_{p\infty}$\
or $L_{p1}$\ and non-com\-mut\-at\-ive generalizations,\sl\ Math.\
Ann.\ {\bf 276} (1986),
105--136.\cr
T & M.~Talagrand,\rm\ The canonical injection from $C([0,1])$\ into
$L_{2,1}$\ is not of cotype~2,\sl\ Contemporary Mathematics, Volume
{\bf
85} (1989), 513--521.\cr
}

\comment{
\vfill

\hfil\vbox{\halign{#\hfil\cr
S.J.~Montgomery-Smith,\cr
Department of Mathematics,\cr
University of Missouri,\cr
Columbia, MO 65211.\cr
\cr\cr\cr\cr\cr\cr
}}
\hfil\vbox{\halign{#\hfil\cr
M.~Talagrand,\cr
Department of Mathematics,\cr
The Ohio State University,\cr
231 W.\ 18th Avenue,\cr
Columbus, OH 43210.\cr
\cr
Equipe d'Analyse -- Tour 46,\cr
Universit\'e Paris VI,\cr
4 Place Jussieu,\cr
75230 Paris Cedex 05.\cr
}}
\hfil

\vskip 0.5 true in

\eject
}

\bye